# DISCRETE EVENT SIMULATION TO EVALUATE SHELTER CAPACITY EXPANSION OPTIONS FOR LGBTQ+ HOMELESS YOUTH

| | |
|---|---|
| Yaren Bilge Kaya | Renata Konrad |
| Sophia Mantell | Andrew C. Trapp |
| Kayse Lee Maass | |
| | |
| Mechanical and Industrial Engineering | School of Business |
| Northeastern University | Worcester Polytechnic Institute |
| 360 Huntington Ave | 100 Institute Rd |
| Boston, MA 02115, USA | Worcester, MA 01609, USA |
| Geri L. Dimas | Meredith Dank |
| | |
| Data Science Program | Marron Institute of Urban Management |
| Worcester Polytechnic Institute | New York University |
| 100 Institute Rd | 60 5$^{th}$ Ave 2$^{nd}$ floor |
| Worcester, MA 01609, USA | New York, NY 10011, USA |

**ABSTRACT**

The New York City (NYC) youth shelter system provides housing, counseling, and other support services to runaway and homeless youth and young adults (RHY). These resources reduce RHY's vulnerability to human trafficking, yet most shelters are unable to meet demand. This paper presents a Discrete Event Simulation (DES) model of a crisis-emergency and drop-in center for LGBTQ+ youth in NYC, which aims to analyze the current operations and test potential capacity expansion interventions. The model uses data from publicly available resources and interviews with service providers and key stakeholders. The simulated shelter has 66 crisis-emergency beds, offers five different support services, and serves on average 1,399 LGBTQ+ RHY per year. The capacity expansion interventions examined in this paper are adding crisis-emergency beds and psychiatric therapists. This application of DES serves as a tool to communicate with policymakers, funders, and service providers—potentially having a strong humanitarian impact.

## 1 INTRODUCTION

Youth homelessness is a severe systemic issue in New York City (NYC), with recent levels being the highest since the great depression (Coalition for the Homeless 2022). According to a 2018 study from Chapin Hill, approximately 4,500 RHY are homeless each night in NYC and there is a severe gap in youth-specific services—such as prevention systems and shelters (Morton et al. 2019). The NYC Department of Youth and Community Development (DYCD) funds eight drop-in centers and 293 crisis beds for youth ages 14 to 24 in NYC; the former of which provides immediate services to RHY, and the latter house youth for up to 120 days (NYC Department of Youth and Community Development, Runaway and Homeless Youth Services 2022). Although services are available for RHY in NYC, there is not nearly enough capacity to meet the demand, and services are often insufficient in matching their needs (Clawson et al. 2009).

Runaway and homeless youth and young adults (RHY) are highly vulnerable to human trafficking—defined as the exploitation of a person, usually in the form of sex or labor, for monetary gain or benefit (Gajic-Veljanoski and Stewart 2007). A large 10-city study, conducted from 2014 to 2016, found that 19% of RHY had been victims of some form of human trafficking and that 68% of youth who had been trafficked or engaged in commercial sex had done so while homeless (Murphy 2016). Another study focusing on youth at the NYC youth shelter Covenant House found that 48% of youth who engaged in commercial sex did it because they did not have a place to stay (Bigelsen and Vuotto 2013). The vulnerability to trafficking increases even more for LGBTQ+ youth; the previously mentioned 10-city study also found that 24% of LGBTQ+ RHY were victims of sex trafficking (Murphy 2016). LGBTQ+ youth account for up to 40% of the RHY population, although only 3-5% of the general youth population identify as LGBTQ+ (Xian et al. 2017). This further highlights the importance of LGBTQ+ specific RHY services and shelters. Additionally, homelessness greatly decreases opportunities for youth seeking employment or education (Bigelson and Vuotto 2013). When youth are unable to find licit employment opportunities, they are much more likely to be trafficked for sex or labor to meet their basic needs (Xian et al. 2017).

Shelter and support services have the potential to help reduce RHY vulnerability to trafficking and exploitation (Ide and Mather 2019). Access to safe, reliable housing is one of the most effective interventions in preventing youth from being exploited and trafficked (Murphy 2016; Dank et al. 2015). Bruhns et al. 2018 conducted a series of interviews in 2018 with survivors of childhood exploitation and human trafficking who highlighted that leaving trafficking requires external assistance from "comprehensive, nonjudgmental services" for the entire duration of their post-trafficking experience. The survivors clearly stated that they needed their basic material and safety needs met (that is, reliable housing and a stable income) to realistically exit the human trafficking situation (Bruhns et al. 2018). Beyond shelter, support services—such as case management and psychiatric care—can further decrease RHY vulnerability to trafficking (Bigelson and Vuotto 2013). Housing and service interventions are both an immediate tool for reducing RHY trafficking vulnerability, as well as a long-term investment in RHY health by decreasing the likelihood of chronic adult homelessness (Hsu et al. 2021). Thus, providing and improving access to support services for RHY, especially those with marginalized identities, is critical.

Given this context, our goal is to model RHY access to shelter resources and analyze potential improvements. This is accomplished via a discrete event simulation (DES) of a shelter in NYC specifically serving LGBTQ+ RHY. The simulation models the process of RHY entering the shelter and using housing and support services based upon their specific needs. The model is built with a combination of secondary data collected from publicly available resources and primary data obtained in interviews with service providers and key stakeholders, such as NYC Coalition for Homeless Youth. It incorporates realistic RHY behaviors such as queue abandonment, varied stay times, and individualized resource needs. This provides a platform for testing capacity expansion interventions and analyzing the potential impact on the shelter and RHY.

The modeled shelter has 66 crisis emergency beds along with five support services: case management, drug counseling, health insurance enrollment, psychiatric services, and medical services. This paper replicates the current state of the simulated shelter and analyzes potential impacts of two interventions: increasing the number of (1) crisis-emergency beds and (2) psychiatric therapists available. This application of discrete event simulation demonstrates expanding shelter services in NYC as an effective way to meet the needs of many RHY, which has a great humanitarian impact.

## 2  RELATED WORKS

Operations research techniques have been applied to several problems regarding access to housing for populations vulnerable to or experiencing homelessness. Kaplan (1987) uses a queueing theory model to analyze equity of wait times for public housing assignments based on the races of tenants. Johnson and Hunter (2000) provide an optimization model to maximize the social benefit and equity of government-issued rental voucher programs within program and policy constraints. Johnson and Smilowitz (2012) create a decision model for finding affordable housing in Pittsburgh, Pennsylvania. The model involves a user interface which allows users to input decisions on certain factors, which improves the autonomy and accessibility for low-income families searching for affordable housing. Johnson (2007) proposes a model to optimize the locations of affordable housing projects by maximizing both the efficiency and equity of the project locations. All of these models analyze various programs that aim to decrease homelessness by expanding accessibility to housing.

Analytical models may be paired with policy simulations to examine potential systematic improvements and increase effectiveness in policy making for homelessness. Early (1999) creates a microeconomic predictive model and a policy simulator to examine the causes of homelessness and best interventions. This analyzed the best availabilities of various assistance programs such as homeless shelters, rent assistance, and inexpensive traditional housing. Azizi et al. (2018) proposes a mixed-integer optimization framework to test the predicted effectiveness of policies for allocating youth housing resources. It successfully determined policy alternatives that increased fairness for RHY of different demographics without reducing efficiency. Kaya et al. (2022) uses a linear programming model to evaluate effective capacity expansion techniques for the NYC system of shelter organizations for RHY. All of these formulations can aid in future policy design and resource allocation decisions.

A handful of existing simulation models address the healthcare needs within homeless populations. Brewer et al. (2001) propose a Markovian simulation of the transmission of tuberculosis within homeless populations. This model examines the impact of several control strategies on cases and deaths within the population, finding that a 10% increase in treatment access had a greater effect on decreasing tuberculosis cases and deaths than any improvements in care quality. Reynolds et al. (2010) provide a discrete event simulation of a healthcare clinic serving homeless patients with a large variety of needs. The clinic specifically provides care for homeless individuals and relies entirely on volunteer staff. The simulation model seeks to improve processes within the clinic to use its limited resources efficiently and divert excess patients from nearby emergency departments. Chapman et al. (2021) conduct a study of outbreaks of COVID-19 in homeless shelters using an individual-state microsimulation. The simulation found low success in any control strategies in preventing outbreaks and emphasizes the need for non-congregate housing in the event of public health emergencies. Homeless populations have unique healthcare needs as opposed to general populations; Reynolds et al. (2010), Brewer et al. (2001), and Chapman et al. (2021) all found validity in using simulation modelling techniques to address health issues within the context of homelessness.

The model proposed in this paper is a unique application of discrete event simulation to analyze the capacity of a singular crisis shelter for LGBTQ+ homeless youth. Azizi et al. (2018) and Kaya et al. (2022) are the only aforementioned papers that focus on homeless youth. The increased vulnerability to human trafficking that homeless youth, especially LGBTQ+ RHY, face poses a dire need for understanding and improving available resources. Although DES models of RHY exist, this method has not used to analyze potential capacity expansion. In contrast, DES has been used extensively to evaluate capacity expansion in healthcare and public programs (Zhu et al. 2012; Ia et al. 2016; Ferreira et al. 2008; Fraga et al. 2016; Park and Noh 1986). Additionally, while much of the literature is situated in the longer-term solution of affordable housing, (Johnson and Smilowitz 2012; Kaplan 1987; Johnson and Hunter 2000; Johnson 2007), there is a gap in capacity expansion of short-term housing programs. This simulation model analyzes crisis bed space availability and other immediate crisis interventions for RHY. To the best of our knowledge, DES has yet to be applied for capacity expansion interventions within a LGBTQ+ youth crisis shelter. DES is an effective tool for testing several potential capacity expansion interventions and communicating the results to policymakers, funders, and service providers of RHY shelter services.

## 3   DATA AND COMMUNITY PARTNERS

This section briefly presents our data acquisition process. The model uses a combination of secondary data collected from publicly available resources and primary data obtained in interviews and surveys with service providers and key stakeholders, such as NYC Coalition for Homeless Youth and NYC Mayor's Office Youth Homelessness Taskforce. The service provider survey yields data regarding shelter services offered, existing capacities, approximate lengths of stay for RHY, and approximate arrival rates.

Information from stakeholder interviews and publicly available data was used to create individual needs-profiles for RHY that we use in our model, which include information about support services needed.

Housing and support services are provided to RHY to ensure ongoing safety, serve basic needs, and reduce vulnerability to exploitative experiences such as sex and labor trafficking. One-third of RHY enter the shelter requesting a crisis bed and there are 66 available beds in the shelter. All RHY in the system seek support services – five support services in addition to housing: case management, drug counseling, health insurance enrollment, psychiatric services, and medical services. Medical services refer to the procedures that are conducted by travelling nurses and doctors in the shelter. These procedures include physical check-ups, testing for sexually transmitted infections and COVID-19, and vaccination updates. All support services, with the exception of case management, are provided by contractors who work at shelters for a limited number of hours per week, thus limiting access to services. The estimated capacity of services within the current system of the shelter is shown in Table 1.

The needs-profile data informs the demand for shelter resources by the homeless youth population. A needs-profile is made up of integers representing the monthly quantity of appointments a youth requests for each service type. These integers are randomly generated based on demand probabilities for each service. For example, a youth's needs-profile may indicate they do not need drug counseling but need to see a psychiatrist weekly. Accordingly, the youth entity will seize zero drug counseling appointments and four psychiatric service appointments per month. We assume that time conflicts do not exist when scheduling these service appointments and youth are able to use the support service appointments that they claim; the model focuses on capacity rather than scheduling. The proportion of youth using each service and the range of monthly appointment quantities per youth are available in Table 1.

Table 1: Estimated Service Capacities.

| Support Service | Total Appointments Available per Month | Percentage of Youth Requesting Service | Appointments per Youth per Month |
|---|---|---|---|
| Case Management | 400 | 100% | 2 - 4 |
| Drug Counselling | 60 | 40% | 0 - 4 |
| Insurance Enrollment | 34 | 50% | 0 - 1 |
| Psychiatric Services | 56 | 50% | 0 - 4 |
| Medical Services | 192 | 90% | 1 - 5 |

## 4 METHODS

### 4.1 Discrete Event Simulation Model

There are two types of entities in our model of the shelter: (i) youth who are seeking to receive a crisis emergency bed *(bed-seeking youth, or BSY)* and (ii) youth who are only seeking support services *(non-bed-seeking youth, or NBSY)*. Figure 1 shows the logic of how those entities flow through the system.

Youth arrive in the shelter with different personal attributes, such as age and the list of services they need. These characteristics determine certain behaviors in the system such as the number of days a youth remains in the shelter and their willingness to wait for services within the system; we call these Stay Attributes, outlined in Table 2. When an entity enters the shelter, we assign its initial attributes. For an NBSY entity, they continue directly to requesting their support service needs. For BSY entities, they first request a bed resource. Once they are able to seize a bed, they move onto requesting their support service needs. If an entity waits in the queue longer than their Bed Patience attribute value, they will renege from the bed queue. While we do not know the exact number, interviews with service providers justified assuming that one quarter of BSY who renege from the bed queue will exit the shelter entirely and three quarters remain in the system to use other support services, thereafter, functioning as NBSY. All entities, regardless of type, are assigned service needs-profiles based on demands discussed in Section 3, Table 1.

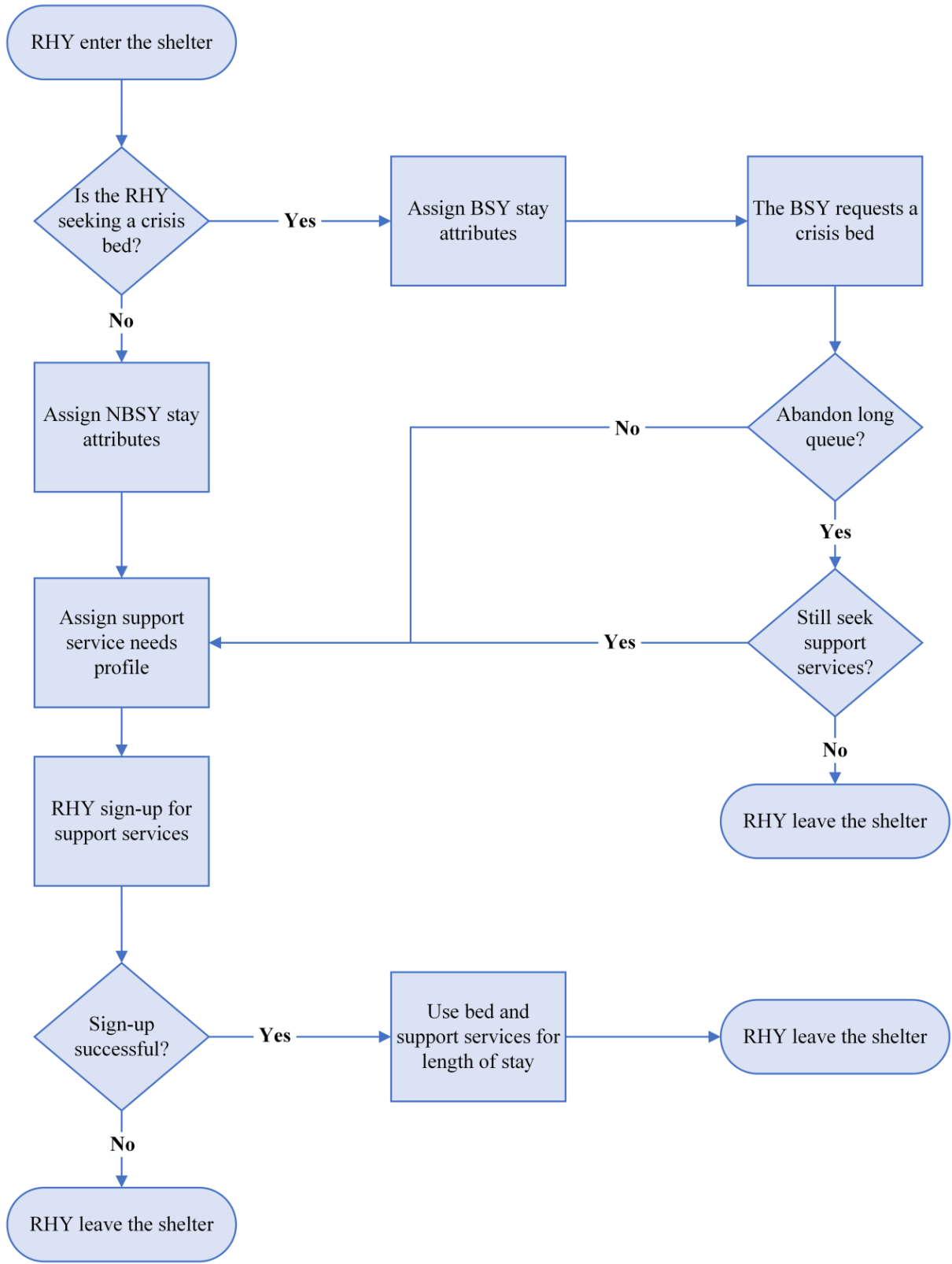

Figure 1: Shelter Intake Process Flow.

Table 2: Youth Stay Attributes – Assigned Values in Triangular Distribution of (days).

| Attribute | Definition | Group | Triangular Distribution (Days) |
|---|---|---|---|
| Length of Stay | The total time a youth stays in the system while using any requested services | BSY, 16-20 y.o.<br>BSY, 21-24 y.o.<br>NBSY | (30, 75, 90)<br>(60, 120, 180)<br>(7, 14, 30) |
| Bed Patience | The maximum time a youth will wait in the crisis bed queue before reneging | BSY, all ages | (3, 5, 7) |
| Service Patience | The maximum time a youth will wait in a support service queue before reneging | All | (1, 7, 14) |

The entity then enters the support service sign-up process of the model. This subsystem operates under two assumptions: (1) if a service does not have any available appointments, the youth will enter a waitlist to receive that service and (2) a youth is able to be on multiple waitlists at a time. To implement the model logic, an original entity is duplicated four times to represent five concurrent service assignments, as seen in Figure 2. Each of the five identical entities undergoes a parallel service sign-up process to acquire appointments for each desired support service simultaneously.

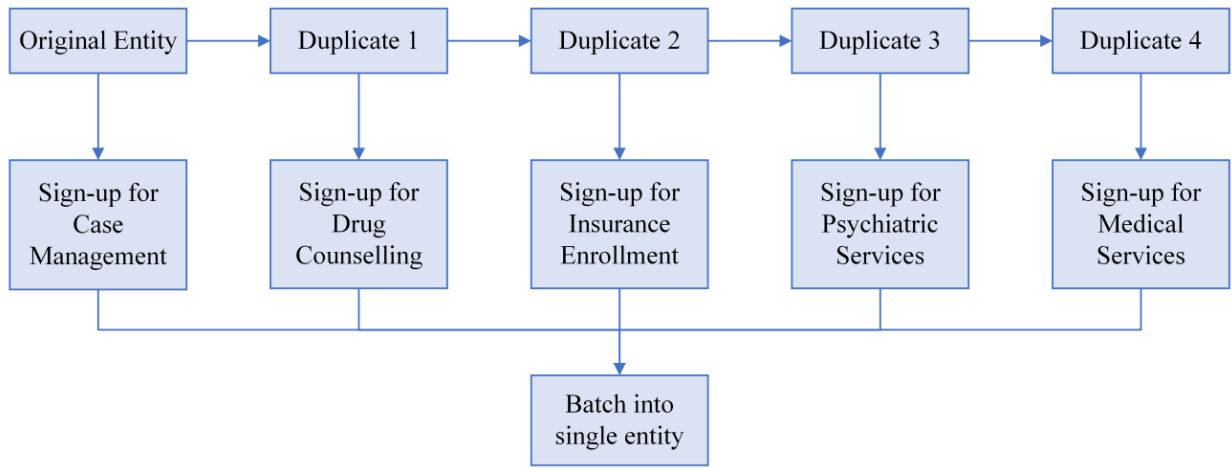

Figure 2: Duplicates Proceed through Support Service Sign-ups in Parallel and then Batch.

The sign-up process for support services is shown in Figure 3 below, with psychiatric services as an example. If an entity's needs-profile does not list a particular service, the entity will bypass the queue and not claim any appointments of that service. Otherwise, the entity will wait in the queue to seize the number of appointments it needs of that service. If the entity waits for any service longer than their *Service Patience* attribute, they will abandon the queue and bypass that service. Once all duplicates have gone through their respective service sign-up process, they are batched back into one entity and sent to the final stage of the model.

Any youth that reneges from all assigned services—and thus has not claimed any appointments—leaves the shelter. Otherwise, the entity stays in the system and uses any bed or support services they have seized for the remainder of their assigned length of stay. Once the entity has stayed the duration of its stay time, it releases all appointments and beds it has claimed and leaves the system. Our model assumes that the system does not track repeated visitors. That is, if a youth leaves the shelter but returns months later to receive more services, the simulation processes them as an entirely new entity.

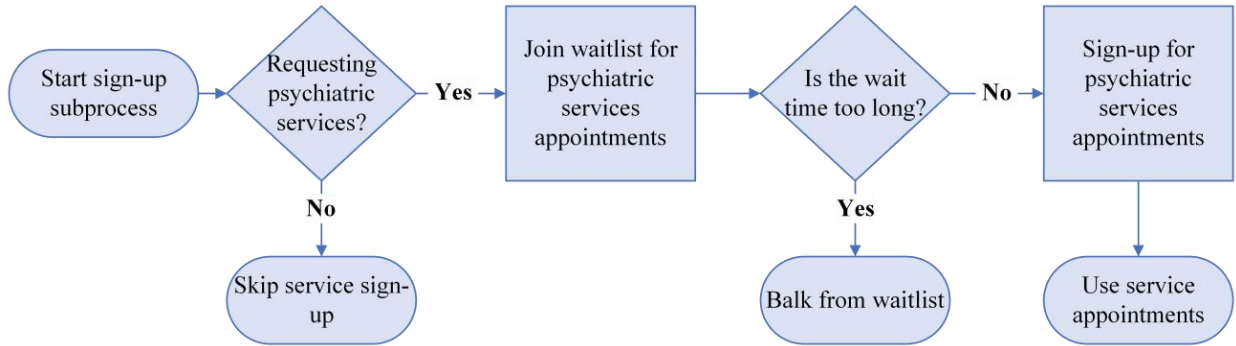

Figure 3: Psychiatric Services Sign-Up Process Flow.

## 4.2 Verification and Validation

To verify the simulation model, we consulted service providers regarding our logical process flow, assumptions, and outputs. Each shelter organization in NYC is unique in the population it serves, services it provides, and its process flow, which creates barriers for the model's transferability. Therefore, several nuances may need to be altered to generalize and apply the model to other shelters in NYC. Nonetheless, service providers verified that inputs and outputs of the simulation do reflect the real-world conditions. We intend to validate specific outputs as more primary data becomes available.

## 5 RESULTS

### 5.1 Baseline Run

We ran the model for a one year warm-up period prior to collecting statistics for another full year run with 100 replications. In that year, an average of 1,399 total youth (469 BSY and 930 NBSY) entered the shelter. The shelter served an average of 1280 youth in that one year period while 119 abandoned the crisis bed queue and left the system without being served. Table 3 below shows the average and maximum wait times of youth who did not abandon the queue, the average utilization, and percentage of RHY who abandoned the queue for each of the shelter's services.

Table 3: Service Wait, Utilization, and Queue Abandonment Statistics in a One Year Period.

| Service | Avg. Wait Time (Days) | Max. Wait Time (Days) | Avg. Utilization | Percent Reneged |
|---|---|---|---|---|
| Crisis Bed | 3.01 | 8.4 | 99.8% | 25.3% |
| Case Management | 1.84 | 9.4 | 98.3% | 54.5% |
| Drug Counseling | 5.89 | 19.6 | 91.1% | 63.8% |
| Insurance Enrollment | 5.4 | 16.4 | 96.7% | 63.6% |
| Psychiatric Services | 8.79 | 24.7 | 92.7% | 68.9% |
| Medical Services | 3.98 | 12.9 | 97.6% | 56.4% |

### 5.2 Crisis Bed Capacity Expansion

Crisis beds are the most critical resource in the current system, as BSY may not have another place to sleep for the night. Of the average 119 youth who reneged from the crisis bed queue in a year, 28 left the shelter entirely while 91 stayed in the system and sought other support services. These results motivated us to explore the impact of interventions to improve access to shelter beds and services for youth. Table 4 and Figure 4 show the potential impact of the shelter acquiring more space and funding to increase the number of beds. Adding five beds to the shelter (that is, increasing from 66 to 71 beds) does not appear to have a strong impact on youth wait times or queue abandonment. However, adding fifteen beds (that is, increasing to 81 beds) reduces the percent of BSY who renege to less than half of its current value and decreases the average wait time for a bed by 31%. Furthermore, adding 25 beds (that is, increasing to 91 beds) reduces the average wait time to under a day (0.8 days) and leaves only 4% of youth reneging due to the long wait time. If the shelter were able to add 40 beds to its existing capacity, the maximum wait time for a bed decreases to 0.1 days, which is less than the patience threshold for any youth, leading to zero youth reneging. This presents an ideal scenario for the supply of beds to match demand. Although space and financial constraints make it unattainable for most shelters to add 40 beds, any additional beds can provide more RHY with safe housing. Adding adequate beds to reduce wait times has the potential to greatly decrease RHY vulnerability to exploitative situations. The greatest need is for crisis beds due to the current high demand, long wait times, and queue abandonment rates, and the magnitude of impact on youth safety provided by accessing stable housing.

Table 4: Crisis Emergency Bed Capacity Expansion Interventions

| Intervention | Quantity of Beds | Avg. Wait Time (Days) | Max. Wait Time (Days) | Avg. Utilization | Percent Reneged |
|---|---|---|---|---|---|
| Baseline | 66 | 3.0 | 8.4 | 100% | 25% |
| Adding 5 beds | 71 | 2.7 | 8.9 | 99% | 24% |
| Adding 10 beds | 76 | 2.4 | 11.0 | 97% | 14% |
| Adding 15 beds | 81 | 2.1 | 8.9 | 98% | 10% |
| Adding 20 beds | 86 | 1.1 | 6.6 | 96% | 4% |
| Adding 25 beds | 91 | 0.8 | 7.5 | 93% | 4% |
| Adding 30 beds | 96 | 0.7 | 6.6 | 86% | 3% |
| Adding 35 beds | 101 | 0.3 | 5.9 | 85% | 1% |
| Adding 40 beds | 106 | 0.0 | 0.1 | 83% | 0% |

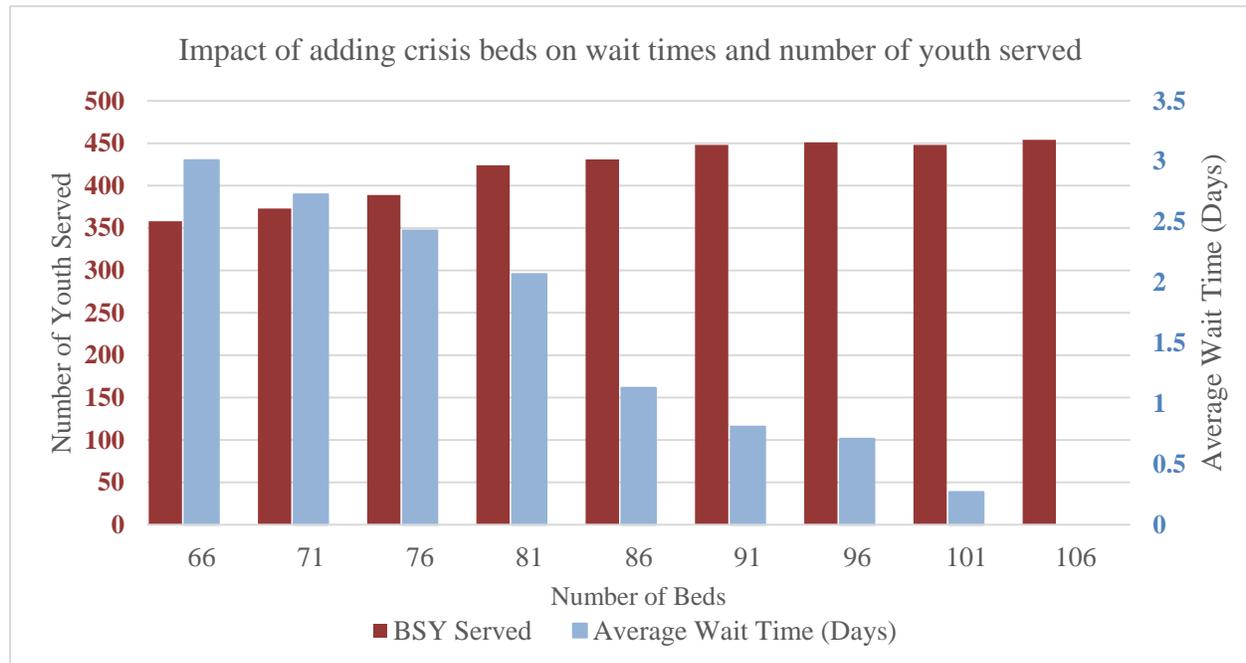

Figure 4: Impact of Crisis Emergency Bed Capacity Expansion Interventions

## 5.3 Psychiatric Services Capacity Expansion

The shelter we are modeling provides mental health support in various forms to RHY such as group therapy, weekly individual therapy, or seeing a psychiatrist. Even though this is not a mandatory service, it plays a significant role in reducing youth's vulnerability to trafficking and exploitation, hence the service request proportion is at 50% in Table 1. With the existing capacity, 69% of youth renege from the psychiatric service sign-up. Thus, we tested interventions for expanding psychiatric services to determine how many therapists the shelter needs to hire so the wait time will be negligibly short. The current capacity of 56 appointments per month is calculated based on having four part time therapists and psychiatrists available four hours per day, two days per week. Table 5 shows the impact of hiring more therapists.

Table 5: Psychiatric Service Capacity Expansion Interventions.

| Intervention | Appointments per Month | Avg. Wait Time (Days) | Max. Wait Time (Days) | Avg. Utilization | Percent Reneged |
|---|---|---|---|---|---|
| Baseline | 56 | 8.8 | 24.7 | 93% | 69% |
| Adding 1 Part Time Therapist | 72 | 5.9 | 16.8 | 94% | 66% |
| Adding 2 Part Time Therapists | 88 | 4.5 | 15.2 | 93% | 62% |
| Adding 3 Part Time Therapists | 104 | 3.7 | 13.8 | 96% | 53% |
| Adding 4 Part Time Therapists | 120 | 2.3 | 11.5 | 94% | 52% |
| Adding 5 Part Time Therapists | 136 | 1.4 | 8.2 | 88% | 50% |
| Adding 6 Part Time Therapists | 152 | 0.5 | 12.1 | 78% | 50% |
| Adding 7 Part Time Therapists | 168 | 0.2 | 4.5 | 73% | 49% |

If the shelter has the resources to hire more therapists, the impact on wait times and queue abandonment rates is noticeable by hiring just four more therapists. By increasing the number of part time therapists from 4 to 10, on average, youth will wait less than a day to be assigned. Considering the importance of mental health support to RHY, it is a worthwhile capacity expansion which decreases the wait time to acquire psychiatric appointments and thus increase youth access.

## 6 CONCLUSIONS AND FUTURE WORK

This paper presents a discrete event simulation model of a crisis shelter for LGBTQ+ runaway and homeless youth and young adults in New York City. The model simulates the process of youth seeking, acquiring, and using the services that they need within the shelter. Services include: a crisis emergency bed, case management, drug counseling, insurance enrollment, psychiatric care, or medical assistance. The simulation is useful for evaluating bottlenecks and shortcomings of the current system and for testing interventions to expand the capacity of certain services. Through expanding the capacity of the services within the shelter, more RHY will be able to access shelter support services with less wait time, reducing their vulnerability to exploitation and improving their quality of life.

Applying DES to youth homelessness services presented several challenges and limitations. While we took care to base model inputs and parameters on realistic data and conversations with youth and service providers, primary data collection was impacted by COVID-19, and therefore as additional data becomes available the model can be updated to reflect the most recent trends. Additionally, while evaluating the capacity expansion interventions, we assume that the demand for crisis-emergency beds and support services remain constant, since increasing the availability of shelter services would not create more homelessness and increase the demand. Hence, exploring the effect capacity expansion has on demand is an area left for future study.

Future work can expand upon the structure of the presented model to a variety of different applications. For example, the model can be adapted and applied to similar crisis housing and drop-in shelters by simply altering input parameters and the services offered. Furthermore, the model could be expanded to simulate multiple shelters simultaneously and provide a representation of the youth shelter system in NYC (or other cities) as a whole. Youth crisis shelters form a complex network in NYC; RHY may use resources from various shelters or move between shelters depending on availability. Modeling the shelter network is an area of further research with great potential. Applying simulation and other analytical tools to homelessness and human trafficking provides a new perspective on solutions that makes the best use of humanitarian aid given the limited resources available. The results of this simulation model and further work may be used to communicate the importance of expanding youth shelter services to policymakers, funders, and service providers.

## ACKNOWLEDGMENTS

This material is based upon work supported by the National Science Foundation under Grant No. CMMI - 1935602. We thank Andrea Hughes, Jamie Powlovich, and Cole Giannone for their efforts in outreach, data collection, and social justice. Additionally, we are grateful to Isabel Cusack, Chia-Hsiang Hsu Tai, John Hooper, Thomas Minutella, and Ana Sosa for their contributed work on this project.

## REFERENCES

Azizi, M. J., P. Vayanos, B. Wilder, E. Rice, and M. Tambe. 2018. "Designing Fair, Efficient, and Interpretable Policies for Prioritizing Homeless Youth for Housing Resources". In *Proceedings of International Conference on the Integration of Constraint Programming, Artificial Intelligence, and Operations Research*, June 26th-29th, 2018, Delft, The Netherlands, 35-51, Springer International Publishing.

Bigelsen, J. and S. Vuotto. 2013. "Homelessness, Survival Sex, and Human Trafficking: As Experienced by the Youth of Covenant House New York". Covenant House, New York City, New York.

Brewer, T. F., S. J. Heymann, S. M. Krumplitsch, M. E. Wilson, G. A. Colditz, and H. V. Fineberg. 2001. "Strategies to Decrease Tuberculosis in US Homeless Populations: A Computer Simulation Model". *Journal of American Medical Association* 286(7):834–842.

Bruhns, M. E., A. del Prado, J. Slezakova, A. J. Lapinski, T. Li, and B. Pizer. 2018. "Survivors' Perspectives on Recovery from Commercial Sexual Exploitation Beginning in Childhood". *Journal of Counseling Psychology* 46(4): 413–455.

Chapman, L. A. C., M. Kushel, S. N. Cox, A. Scarborough, C. Cawley, T. Q. Nguyen, I. Rodriguez-Barraquer, B. Greenhouse, E. Imbert, and N. C. Lo. 2021. "Comparison of Infection Control Strategies to Reduce COVID-19 Outbreaks in Homeless Shelters in the United States: A Simulation Study". *BMC Medicine* 19(1): 116–116.

Clawson, H., N. Dutch, A. Solomon, and L. G. Grace. 2009. "Human Trafficking into and within the United States: A Review of the Literature". Office of the Assistant Secretary for Planning and Evaluation, U.S. Department of Health and Human Services, Washington, D.C.

Coalition for the Homeless. 2022. "Basic Facts About Homelessness: New York City." https://www.coalitionforthehomeless.org/basic-facts-about-homelessness-new-york-city, accessed 25th February 2022.

Dank, M., J. Yahner, K. Madden, I. Banuelos, L. Yu, A. Ritchie, M. Mora, and B. Connor. 2015. "Surviving the streets of New York Experiences of LGBTQ Youth, YMSM, and YWSW Engaged in Survival Sex." Urban Institute, New York City, New York.

Early, D. W. 1999. "A Microeconomic Analysis of Homelessness: An Empirical Investigation Using Choice-Based Sampling". *Journal of Housing Economics* 8(4):312–327.

Ferreira, R. B., F. C. Coelli, W. C. A. Pereira, and R. M. V. R. Almeida. 2008. "Optimizing Patient Flow in a Large Hospital Surgical Centre by Means of Discrete-Event Simulation Models." *Journal of Evaluation in Clinical Practice* 14(6):1031-1037.

Fraga, C. C. S., J. Medellín-Azuara, and G. F. Marques. 2016. "Planning for Infrastructure Capacity Expansion of Urban Water Supply Portfolios with an Integrated Simulation-Optimization Approach." *Sustainable Cities and Society* 29:247-256.

Gajic-Veljanoski, O. and D. E. Stewart. 2007. "Women Trafficked into Prostitution: Determinants, Human Rights and Health Needs." *Transcultural Psychiatry* 44(3):338–58.

Hsu, H.-T., C. Hill, M. Holguin, L. Petry, D. McElfresh, P. Vayanos, M. Morton, and E. Rice, E. 2021. "Correlates of Housing Sustainability Among Youth Placed into Permanent Supportive Housing and Rapid Re-Housing: A Survival Analysis". *Journal of Adolescent Health* 69(4):629–635.

Ia, E. M., K. E. H. Lich, R. Wells, A. R. Ellis, M. S. Swartz, R. Zhu, and J. P. Morrissey. 2016. "Increasing Access to State Psychiatric Hospital Beds: Exploring Supply-side Solutions." *Psychiatric Services* 67:523–528.

Ide, M. and D. M. Mather. 2019. "The Structure and Practice of Residential Facilities Treating Sex Trafficking Victims". *Journal of Human Trafficking* 5(2):151–164.

Johnson, M. P. and K. Smilowitz. 2012. "*Community-Based Operations Research: Introduction, Theory, and Applications*". In *Community-Based Operations Research*, edited by M. P. Johnson, 3–36. Catonsville, Maryland. The Institute for Operations Research and the Management Sciences.

Johnson, M. P. and A. P. Hurter. 2000. "Decision Support for a Housing Relocation Program Using a Multi-Objective Optimization Model". *Management Science* 46(12):1569–1584.

Johnson, M. P., 2007. "Planning Models for Affordable Housing Development". *Environment and Planning B: Planning and Design* 34(1): 501–523.

Kaplan, E. H., 1987. "Analyzing Tenant Assignment Policies". *Management Science* 33(3):395–408.

Kaya, Y. B., K. L. Maass, G. L. Dimas, R. Konrad, A. C. Trapp, and M. Dank. 2022. "Improving Access to Housing and Supportive Services for Runaway and Homeless Youth: Reducing Vulnerability to Human Trafficking in New York City". https://arxiv.org/abs/2202.00138, accessed 10th April 2022.

Morton, M. H., M. A. Kull, R. Chavez, A. J. Chrisler, E. Carreon, and J. Bishop. 2019. "New York City Youth Homelessness System Map & Capacity Overview". Chapin Hill, University of Chicago, Chicago, Illinois.

Murphy, L. 2016. "Labor and Sex Trafficking Among Homeless Youth: A Ten-City Study". Loyola University, New Orleans, Louisiana.

NYC Department of Youth and Community Development, Runaway and Homeless Youth Services. 2022. "Local Law 86 Report to the Speaker of the City Council. Fiscal Year 2021". New York City, New York.

Park, C. S. and Y. D. Noh. 1986. "An Interactive Port Capacity Expansion Simulation Model". *Engineering Costs and Production Economics* 11(1):109-124.

Reynolds, J., Z. Zeng, J. Li, and S. Chiang. 2010. "Design and Analysis of a Health Care Clinic for Homeless People Using Simulations". *International Journal of Health Care Quality Assurance* 23(6):607–620.

Xian, K., S. Chock, and D. Dwiggins. 2017. "LGBTQ Youth and Vulnerability to Sex Trafficking". In *Human Trafficking Is a Public Health Issue: A Paradigm Expansion in the United States,* edited by M. Chisolm-Straker and H. Stoklosa, 141–152. New York City, New York. Springer International Publishing.

Zhu, Z., B. H. Hen, and K. L. Teow. 2012. "Estimating ICU Bed Capacity Using Discrete Event Simulation". *International Journal of Health Care Quality Assurance* 25(2):134–144.

## AUTHOR BIOGRAPHIES

**YAREN BILGE KAYA** is a PhD Candidate in the Mechancial and Industrial Engineering department at Northeastern University. Prior to joining the PhD program at Northeastern, she received her MSc in Industrial Engineering from the University of South Florida and BSc in Industrial Engineering from Ozyegin University in Istanbul. Improving access to public services and increasing efficiency in healthcare systems have been the main themes of her research. Her email address is kaya.y@northeastern.edu and her website is https://yarenbilgekaya.wixsite.com/yarenbilgekaya.

**SOPHIA MANTELL** is an graduate Civil Engineering, Transportation student at Northeastern University. She completed her undergraduate degree in Industrial Engineering in May 2022 at Northeastern University. She is a student researcher in the Operations Research and Social Justice Lab. She is interested in operations research applications in public programs, social justice, and sustainable transportation. Her email is mantell.s@northeastern.edu.

**KAYSE LEE MAASS** is an Assistant Professor in the Mechanical and Industrial Engineering department at Northeastern University where she leads the Operations Research and Social Justice Lab. She earned her Ph.D. in Industrial and Operations Engineering at the University of Michigan and completed her postdoctoral studies in the Department of Health Sciences Research at the Mayo Clinic. Dr. Maass's research focuses on the application of operations research methodology to social justice, access, and equity issues within human trafficking, mental health, housing, and food justice contexts. Her email address is k.maass@northeastern.edu. Her website is https://kayse-lee.wixsite.com/kayseleemaass.

**RENATA KONRAD** is an Associate Professor of Operations and Industrial Engineering at Worcester Polytechnic Institute. Her research centers on applications of operations research to anti-human traffickign operations and healthcare access. She was awarded the 2021 IISE Modeling and Simulation Division Teaching Award for her innovations in teaching simulation. She earned her PhD in Industrial Engineering from Purdue University. Her email address is rkonrad@wpi.edu and website is https://www.renata-konrad.com.

**ANDREW C. TRAPP** is an Associate Professor of Operations and Industrial Engineering at Worcester Polytechnic Institute, with joint appointments in Mathematical Sciences and Data Science. He researches the use of prescriptive and predictive analytics, together with algorithms, to effectively allocate scarce resources for systems that serve vulnerable peoples. He creates novel analytical technologies and open-source software to improve quality of life, increase fairness, restore dignity, and generate significant societal impact. His email address is atrapp@wpi.edu and his website is http://users.wpi.edu/~atrapp/.

**GERI L. DIMAS** is an Ph.D. Candidate in the Data Science Program at Worcester Polytechnic Institute. She received her M.S. in Applied Statistics from Bowling Green State University and received both a B.A. in Actuarial Science and a B.S. in Computer Science from Roosevelt University. Geri's research focuses on applications of data science at the intersection of societal issues such as immigration, human trafficking and homelessness. Her email is gldimas@wpi.edu.

**MEREDITH DANK** is a Research Professor and directs the Human Exploitation and Resilience initiative of the NYU Marron Institute of Urban Management. She is a nationally recognized expert on human trafficking. She has served as principal investigator for nearly two dozen human trafficking studies, funded by both federal agencies (including the Department of Justice, US Department of State, Department of Health, and Human Services, among others) and private foundations. Her email is mld236@nyu.edu.